\documentclass[11pt]{amsart}
\usepackage[utf8]{inputenc}
\usepackage{ragged2e}
\usepackage{mathtools}
\usepackage{mathrsfs}
\usepackage{float}
\usepackage{tikz}
\usepackage{amsfonts}
\usepackage{amsmath}
\usepackage{enumerate}
\usepackage{amsthm}
\usepackage{amssymb}
\usepackage{dutchcal}
\usepackage{wasysym}
\usepackage[foot]{amsaddr}
\usepackage{bbm}
\usepackage[left=2.5cm,right=2.5cm,top=3.2cm,bottom=3cm]{geometry}
\usepackage{caption}
\usepackage[OT2,T1]{fontenc}
\DeclareSymbolFont{cyrletters}{OT2}{wncyr}{m}{n}
\DeclareMathSymbol{\Sha}{\mathalpha}{cyrletters}{"58}

\DeclareMathOperator{\avg}{avg}
\DeclareMathOperator{\sqf}{sqf}

\usepackage{pdfpages}
\usepackage[natbib, maxbibnames=6, maxcitenames=6]{biblatex}
\makeatletter
\renewcommand{\fnum@figure}{Figure \thefigure}

\newtheorem{teo}{Theorem}
\newtheorem{lem}{Lemma}[section]
\newtheorem{prop}[lem]{Proposition}
\newtheorem{thm}[lem]{Theorem}
\newtheorem{cor}{Corollary}[lem]
\newtheorem{defi}[lem]{Definition}

\theoremstyle{remark}

\newtheorem{reml}[cor]{Remark}

\newtheorem{q}{Question}

\def\R{\mathbb{R}}

\def\Q{\mathbb{Q}}
\def\Z{\mathbb{Z}}

\def\D{\Delta}
\def\F{\mathbb{F}}

\def\h{\mathbb{H}}

\def\p{\mathcal{P}}
\def\f{\mathcal{F}}
\def\g{\mathcal{G}}

\definecolor{linkgreen}{RGB}{45,145,85}
\usepackage[
    colorlinks=true,
    linkcolor=linkgreen,
    citecolor=linkgreen,
    urlcolor=linkgreen
]{hyperref}

\usepackage{cleveref}
\crefformat{section}{\S#2#1#3}
\crefformat{subsection}{\S#2#1#3}
\title{On elliptic curves with prescribed analytic rank at least $2$}
\begin{document}
\author{Francesco Naccarato}
\address{Department of Mathematics, ETH Zürich, Zürich CH-8092, Switzerland.}
\email{francesco.naccarato@math.ethz.ch}
\maketitle

\begin{abstract}
    Zywina recently showed that infinitely many elliptic curves over $\Q$ have rank $2$, settling a long-standing open question. We explore the analytic version of this problem, showing that every elliptic curve over $\Q$ has infinitely many quadratic twists with analytic rank either $2$ or $4$, subject to the Riemann Hypothesis for the $L$-functions of all of its twists. We also obtain, under the same assumption, infinitely many congruent number curves with odd analytic rank between $3$ and $17$.
\end{abstract}

\section{Elliptic curve ranks: a brief overview}\label{1}

Let $k$ be a number field. Recall that the group of $k$-rational points of an elliptic curve $E$ defined over $k$, known as its Mordell--Weil group, is finitely generated. This was shown in 1922 for $k=\Q$ by Mordell \cite{mord}, and extended a few years later to arbitrary number fields and abelian varieties by Weil \cite{weil-mw}. The algebraic rank $r_k(E)=\mathrm{rk}_\Z(E(k))$ is the number of independent generators of the Mordell--Weil group $E(k)$; in other words, \begin{align*}
    E(k)\simeq\Z^{r_k(E)}\oplus E(k)_{\textit{tors}},
\end{align*}where $E(k)_{\textit{tors}}$ is the torsion part. As our focus will be on elliptic curves over $\Q$, we will omit all dependencies on the number field from our notation.

The Birch--Swinnerton-Dyer conjecture predicts that $r(E)$ should be equal to the analytic rank $r_{\text{an}}(E)$, which is the order of vanishing at $s=1$ of the $L$-function $L(E,s)$ of the curve. The latter is the analytic continuation of the Euler product
\begin{align}\label{eq:functional equation}
    \begin{split}
    \prod_{p\nmid N_E}\frac{1}{1-a_p(E)p^{-s}+p^{1-2s}}
     \prod_{p\mid N_E}\frac{1}{1-a_p(E)p^{-s}},\qquad \Re(s)>\frac32,
    \end{split}
\end{align}
where $a_p(E)=p+1-\#E(\F_p)$ at primes of good reduction, and $a_p(E)\in\{0,\pm1\}$ at primes of bad reduction.
Recall that elliptic curves over $\Q$ are \textit{modular} by the celebrated works of Wiles, Taylor--Wiles, and Breuil--Conrad--Diamond--Taylor \cite{wiles,taylor-wiles,bcdt}. This is what allows us to obtain the entire continuation of $L(E,s)$ which makes the analytic rank well defined. Indeed, modularity provides a weight-$2$ cusp form $f=f_E$ such that $L(E,s)=L(f,s)$. The quantity $N_E$ in \eqref{eq:functional equation} is the conductor of $E$, which is the smallest integer $N$ such that there exists such a cusp form of level $N$. 

\subsection{The analytic rank}\label{1.1}

Current methods give us very limited control over the analytic rank---with the exception of its parity, encoded in the root number $\omega(E)=(-1)^{r_{an}(E)}$ (which equals to the eigenvalue of $f_E$ with respect to the Fricke involution of $S_2(\Gamma_0(N_E))$, with its sign flipped \cite[p. 204]{diam}). Indeed, we cannot currently prove that a single elliptic curve $E/\Q$ has analytic rank $r$ for any $r>3$. The main issue is with the lower bound, where numerical methods, such as computing the convergent sums \cite[(2-9)]{dokchitser-values} for the central values of $L^{(n)}(E,s)$ to arbitrary precision, will never suffice in showing its vanishing. Even proving that some elliptic curves have analytic rank $3$---a feat due to Gross and Zagier \cite{gz}---had deep consequences: indeed, Goldfeld \cite{go} had already shown that the existence of such a curve would yield an explicit answer to Gauss's class number problem. 

In terms of upper bound, the main tool we have is a version of the Guinand--Weil \textit{explicit formula} which requires the assumption of the Riemann Hypothesis for $L(E,s)$. Even this method is ineffective, in the sense that the prime cutoff one has to reach in \eqref{eq:expin} to approximate the analytic rank within a given margin depends on the height of the smallest low-lying zero of $L(E,s)$, for which we have no effective lower bound.

Observe that difficulties in computing the analytic rank only arise when it is at least $2$, as otherwise computing the root number $\omega(E)$ suffices in complementing a finite nonvanishing computation. Analytic ranks $0$ and $1$ can indeed be shown to occur infinitely often, even in a given quadratic twist family. This follows from the first-moment asymptotics for central values of Luo--Ramakrishnan \cite[Proposition~3.6]{LR97} and Murty--Murty \cite[Theorem~2]{MM91}, averaged over suitable sets of quadratic twists.

\subsection{The algebraic rank}\label{1.2}

The situation for the algebraic rank is more encouraging, as in order to prove a lower bound one needs to ``simply'' exhibit a set of independent nontorsion rational points. Nonetheless, known algorithms for computing $r(E)$ can be shown to run in finite time only under the assumption of the finiteness of some $p$-part of the Tate--Shafarevich group $\Sha(E/\Q)$. 

It is not too hard to prove that there are infinitely many elliptic curves $E/\Q$ with $r(E)=0$, and the same has been known for rank $1$ since the work of Heegner \cite{heegner}. The case $r(E)=2$ was known, until recently, only under the assumption of the parity conjecture \cite{bj}. In 2025, Zywina \cite{zw} settled it by exploiting recent breakthroughs in additive combinatorics. He later extended his result to elliptic curves with rank exactly $3$ and $4$ \cite{zw2026}. Notice that, if we relax the requirement to $r(E)\ge r$, we obtain a far easier question, since we know how to construct elliptic families $E_t/\Q(t)$ of generic rank up to $18$ \cite{elkies}. The following well-known result \cite[Theorem C]{silv-specialization} then ensures that the algebraic rank almost never decreases under specialization:

\begin{thm}[Silverman specialization theorem]
    Let $E_t/\Q(t)$ be a nonconstant elliptic family of generic rank $r$. Then, for all but finitely many values of $t_0\in\Q$, $E_{t_0}$ is an elliptic curve with $r(E_{t_0})\ge r.$
\end{thm} 

More difficult for the algebraic rank is showing an upper bound. This is usually done by some form of arithmetic descent. Making this technique effective at infinitely many rational fibers of an elliptic family was the main breakthrough in Zywina's works.

\subsection{Main results}

In virtue of the above, it natural to ask the following:
\begin{q}\label{q1}
    For which $r>1$ do there exist infinitely many elliptic curves of (analytic) rank $r$?
\end{q}

The main goal of this work is to explore the analytic variant of Question \ref{q1}: assuming GRH, we establish the existence of infinitely many elliptic curves with $r_{\text{an}}(E)\in\{2,4\}$ in any quadratic twist family over $\Q$.

\begin{teo}\label{mainthm}
    Let $E/\Q$ be an elliptic curve and let $\f_E$ be the relative quadratic twist family. Under the assumption of the Riemann Hypothesis for the $L$-functions of the curves in $\f_E$, infinitely many of them have analytic rank either $2$ or $4$.
\end{teo}

As one might expect, the restriction to even values of the analytic rank in Theorem \ref{mainthm} comes from the fact that, in order to exclude curves of rank $1$, it is convenient (at least for us) to work with twists having root number $+1$. In the case of root number $-1$, we are able to prove a somewhat weaker result for the congruent number family:

\begin{teo}\label{odd}
    Let $E^1: y^2=x^3-x$ be the \textit{congruent number curve.} Under the assumption of the Riemann Hypothesis for the $L$-functions of its twists $E^d$, infinitely many of them have odd analytic rank $3\le r_{\text{an}}(E^d)\le17$.
\end{teo}

\subsection{Notation and conventions}

For a nonzero integer $d$, $\chi_d$ denotes the primitive Dirichlet character associated with $\Q(\sqrt d)$. \\
For $n\in\Z\setminus\{0\}$, $\sqf(n)$ denotes its squarefree part (with sign).\\
For any integrable function $W$ on $\R^N$, its Fourier transform is
\[
\widehat W(\xi)=\int_{\R^N} W(x)\exp\bigl(-2\pi i x\cdot\xi\bigr)\,dx.
\]
For any positive integer $m$, we write $e_m(x)=\exp(2\pi i x/m)$.

\section{Ranks in quadratic twist families}\label{2}

\subsection{Basics of quadratic twists}

Recall that any elliptic curve $E/\Q$ admits a short Weierstrass (affine) model \begin{align}
E: y^2=f(x)\label{eq:ec}\end{align}with $f(x)=x^3+ax+b\in\Z[x]$ and $0\neq\D_E=-16(4a^3+27b^2)$. 
\begin{defi}
    Given $E$ as in \eqref{eq:ec} and $d\in\Z\setminus\{0\}$, the \textit{quadratic twist} of $E$ by $d$ is the elliptic curve with Weierstrass model $$E^d: y^2=x^3+ad^2x+bd^3.$$ 
\end{defi}
Recall that $E^d$ also has the affine model \begin{align}
   dy^2=x^3+ax+b.\label{non-w-model} 
\end{align}With the latter, sending $y\mapsto y\sqrt{d}$ gives an isomorphism with $E$ defined over $\Q(\sqrt{d})$, and in particular an injection \begin{align}E^d(\Q)\hookrightarrow E(\Q(\sqrt d)).\label{twist}\end{align} 
It follows that for any $d,n\in\Z\setminus\{0\}$ we have $E^{dn^2}\simeq_\Q E^d$, so we can define:
\begin{defi}
    Given an elliptic curve $E/\Q,$ its quadratic twist family $\f_E$ is the collection of its quadratic twists by squarefree integers $d$.
\end{defi}

\subsection{Known results}

In virtue of Question \ref{q1}, it is natural to ask how ranks behave in quadratic twist families. On the algebraic side, plenty of work has been done for the $2$-Selmer rank---which is an upper bound for the algebraic rank---until the recent breakthroughs of Smith giving full distribution for most such families \cite{smi}. A corollary of Smith's work is that, for such families, at least $50\%$ of twists (when ordered by modulus of the twisting parameter $d$) have rank $0$ and $100\%$ have rank at most $1$---a result very close to the minimalist conjecture for the family. 

The situation for higher ranks is less clear; the existence of infinitely many twists of rank at least $2$ is known since the work of Mestre \cite{mes}, while infinite subfamilies of rank at least $3$ and even $4$ have only been constructed in some cases \cite{rubin-silv}. In terms of upper bounds, there are heuristics \cite{mppv} suggesting that ranks should be bounded, and---at least for all but finitely many twists---by quite a low cap (around $8$, but depending on the $2$-torsion of the family) \cite{granville+}. The consensus on this, and especially on uniform boundedness for all elliptic curves over $\Q$, appears to be quite split, with evidence going both ways. While the absolute record increased from $28$ to $32$ in the span of two years, the highest known rank for congruent number curves has stood at $7$ for more than $20$ years \cite{rogers}, despite new computational efforts scanning up to much larger twist height \cite{icarm}.

On the analytic side, what we can say in terms of lower bounds comes mostly from a combination of parity with deep results of Gross--Zagier \cite{gz} and Kolyvagin \cite{kol}. In terms of upper bounds, the main tool we have is the Guinand--Weil explicit formula, which requires one to assume the Riemann Hypothesis for the $L$-function of the curve. We will give each of these sets of techniques a more in-depth treatment later in the note.  

The explicit formula is especially powerful when applied to families rather than single curves, and is therefore a crucial tool for \textit{averages,} such as Heath-Brown's result that the average analytic rank of all elliptic curves is at most $2$, then improved to $\frac{25}{14}$ by Young \cite{young}. For quadratic twist families, Goldfeld obtained the bound of $3.25$ \cite{gold}, which Heath-Brown improved to $1.5$ \cite{hb}. Let us remark that the analytic analog of the minimalist conjecture for quadratic twist families over $\Q$ is known as Goldfeld's conjecture, which Smith's results show being implied by the Birch--Swinnerton-Dyer conjecture \cite{smi2}.



\subsection{A subfamily of twists with positive rank}
We now describe how to obtain, inside any quadratic twist family, an explicit infinite subfamily of twists with analytic rank at least $1$. This will then be pushed to $2$ in a suitable subfamily by controlling the root number. This argument is, to the author's best knowledge, due to Le Boudec \cite[p. 7]{leb:avgrk}.
Let $E$ be as in \eqref{eq:ec}, let $N_E$ be its conductor, and let $F(u,v)=u(v^3+au^2v+bu^3)$ be the degree-$4$ homogeneous polynomial such that $F(1,v)=f(v)$. Notice that $E^{F(u,v)}$, taken with the model \eqref{non-w-model}, has the rational point $P=(u^{-1}v,u^{-2})$: this is equivalent to $$F(u,v)=u^4F(1,u^{-1}v)=uv^3+au^3v+bu^4,$$ which holds.
The key observation is then the following:
\begin{lem}\label{fam}
There exists a constant $C_E>0$ depending only on $E$ such that for all pairs of coprime integers $(u,v)$ with $\log\max\{|u|,|v|\}>C_E$, the point $P=(u^{-1}v,u^{-2})\in E^{F(u,v)}(\Q)$ has infinite order.
\end{lem}
\begin{proof} Recall that a rational point on an elliptic curve is torsion if and only if its Néron--Tate height is 0. We now sketch Le Boudec's argument \cite{leb:avgrk} for the reader's convenience: \cite[Lemma 3]{leb:ht} allows us to bound the absolute difference between the naive and the Néron--Tate heights of a rational point on any quadratic twist of $E$ by a constant $C_E$ depending only on $E$; we then obtain, by definition of the Weil height, that $P$ is nontorsion for $\log\max\{|u|,|v|\}>C_E$.\end{proof}
\begin{reml}\label{pos}
    The proof shows that the result holds even if $u$ and $v$ are not coprime, as long as $\log\max\{|u|,|v|\}>C_E+\log\gcd(u,v)$.
\end{reml}
\begin{cor}\label{posanal}
    For $C_E$ as in Lemma \ref{fam} and $u,v$ as in Remark \ref{pos}, we have $r_{\text{an}}(E^{F(u,v)})>0$.
\end{cor}
\begin{proof}
    By the works of Gross--Zagier and Kolyvagin \cite{gz,kol}, the Birch--Swinnerton-Dyer conjecture is known for elliptic curves over $\Q$ with analytic rank smaller than 2. So, in our case, the analytic rank cannot be $0$.
\end{proof}

\section{Overview of the method and one prototype case}

\subsection{Averages}

Suppose we are given a family $\f$ of rational elliptic curves and a surjective function $i:I\to\f$, with $I\subset\Z^N$ a subset of positive lower density, for some $N\ge1$ (for us, $N$ will be either $1$ or $2$). We will use the following notion of \textit{average:}

\begin{defi}\label{def:avg}
    Let $W\in C^1_c(\R^N)$ be nonnegative and such that $\sum_{x\in I\cap[T]}W(x/T)\gg T^N$, where $[T]$ denotes the set $[-T,T]^N\cap\Z^N.$ Given a function $R:\f\to\R,$ its $W$-average is defined as \begin{align}
        \avg_{W}R:=\limsup_{T\to\infty}\frac{\sum_{x\in I\cap[T]}R(i(x))W(x/T)}{\sum_{x\in I\cap[T]}W(x/T)}.
    \label{eq:avgdef}\end{align}
\end{defi}

We call $i$ \textit{good} if each one of its fibers has density $0$ in $\Z^N$. The following result follows almost by definition:

\begin{lem}
    Assume $i$ is good. If the inequality $\avg_W R<A$ holds for some $A>0$, there exist infinitely many curves $E\in\f$ with $R(E)<A$.\label{lem:ineq}
\end{lem}
\begin{proof}
    Let $\f'\subset\f$ be the subfamily where $R<A$, and let $I'=i^{-1}(\f')$. If $\f'$ were finite, then $i$ being good would imply $\#(I'\cap[T])=o(T^N)$. As $W$ is bounded and $R$ takes only finitely many values on $i(I')=\f'$, the contribution of $I'$ to the sums at the numerator and denominator of \eqref{eq:avgdef} would also be $o(T^N)$, allowing us to disregard them in the limit and compute both sums for $I\setminus I'$. Since here $R(i(x))\ge A$, this contradicts the hypothesis that $\mathrm{avg}_WR<A.$
\end{proof}

\subsection{The explicit formula}\label{section:expfor}

Let us recall the Guinand--Weil \textit{explicit formula} (see \cite{guin}, \cite{weilexp}), which gives an upper bound for the analytic rank under the assumption of the Riemann Hypothesis for $L(E,s)$. Let \begin{align*}
        U_k(X) = U_k(E,X) = -\sum_{p^k\le X, p\ge5}\frac{a_{p^k}}{p^k}h_X(k\log p)\log p,  \end{align*}
where here $a_{p^k}(E)$ denotes the trace of the $k$-th power of Frobenius; in particular, $a_{p^2}(E)=a_p(E)^2-2p$ at primes of good reduction.
\begin{prop}
Suppose that all nontrivial zeroes $\rho$ of $L(E,s)$ have real part $1$. Let $h:\R\longrightarrow\R$ be even, piecewise $C^1$ and have compact support $K\subset[-1,1]$ and modulus bounded by $1$. Let $\hat{h}$ be its Fourier transform and set $h_X(t)=h(t/\log X)$. Then, for any $X\ge2$: \begin{align}\begin{split}&r_{an}(E)\hat{h}(0) + \sum_{\tau=\Im\rho\neq0}\hat{h}\left(\frac{\tau\log X}{2\pi}\right) =\notag\end{split}\\ =h(0)\frac{\log N_E}{\log X}& + \frac{2}{\log X}\left(U_1(E,X)+U_2(E,X)\right) + O\left(\frac{1}{\log X}\right).\label{eq:expfor}\end{align}
\end{prop}
\begin{cor}
Suppose $L(E,s)$ satisfies the Riemann Hypothesis. Then, for $X\ge2$, we have: \begin{align}r_{an} < \frac{\log N_E}{\log X} + \frac{2}{\log X}\left(U_1(X)+U_2(X)\right) + O\left(\frac{1}{\log X}\right).\label{eq:expin}\end{align}
\end{cor}
\begin{proof} Setting $h(t) = \max(0, 1-|t|)$ we get $\hat{h}(t)=\left(\frac{\sin \pi t}{\pi t}\right)^2$, which is nonnegative. Since both $h$ and $\hat{h}$ evaluate to $1$ for $t=0$, we get \eqref{eq:expin} by neglecting the summation on the left-hand side of \eqref{eq:expfor}. \end{proof}

\begin{reml}
    Applying this formula to a \textit{single} curve yields $$r_{an}(E)\ll\frac{\log N_E}{\log\log N_E},$$see \cite{brum}.
\end{reml}
 For quadratic twist families, we can get cancellation inside of $\sum U_i(X)$ thanks to the following well-known transformation law \cite{hb}:
\begin{lem}
    Let $d$ be a squarefree integer. For $p$ coprime to $2N_E$, we have
\begin{align}
    a_p(E^{d})=a_p(E)\genfrac(){}{0}{d}{p}.
\end{align}\label{lem:apa1}
\end{lem}
The following result, due to Goldfeld \cite{gold}, will allow us to forget about the $U_2$ term in \eqref{eq:expfor} when averaging over the family: 
\begin{lem}Let $d$ be a squarefree integer. For any $X\ge2$, we have
    \begin{align*}
        U_2(E^{d},X)=\frac{1}{4}\log X + O_E(\log\log |d|).
    \end{align*}\label{lem:gold}
\end{lem}
Controlling $U_1$ is harder: using the Hasse bound and the Prime Number Theorem we get \begin{align}
   |U_1(X)|\le 2\sum_{2\le p\le X} \frac{\log p}{\sqrt{p}} \ll \sqrt{X},\label{eq:u1}
\end{align}
which is not particularly useful since its denominator in the explicit formula is $\log X$. In what follows, we will see how to exploit Lemma \ref{lem:apa1} to obtain the cancellation necessary to prove a suitable average-rank statement.
\subsection{The case of the congruent number family}\label{cn}

Among quadratic twist families, one of the more well-known ones is that of \textit{congruent number curves:} quadratic twists of the curve $E:y^2=x^3-x$. The name comes from the 3000-year-old congruent number problem, which asks for an algorithm to determine whether a given positive integer $n$ is the area of a right triangle with rational sidelengths. It turns out that this condition is equivalent to $E^n$ having positive rank. In 1983, Tunnell \cite{tunn} gave a partly conditional solution to this problem, building on Shimura and Waldspurger's theory of half-integer weight modular forms and their connection to special values of $L$-functions.

The congruent number curves have complex multiplication by $k=\Q(i),$ so their $L$-functions are $L$-functions for Hecke characters of $k$. This allows for more direct control on the sums appearing in \eqref{eq:expfor}. In this section we will give a relatively simple proof---initially presented in the author's Master's thesis \cite{nac-thesis} with a slightly weaker bound---that infinitely many congruent number curves have analytic rank in $\{2,4,6\}$, conditional on the Riemann Hypothesis for their $L$-functions. In the next section we will extend and strengthen the result by invoking some tools from Fourier analysis.

\begin{prop}\label{thmcn}
    Assuming the Riemann Hypothesis for the $L$-functions of the congruent number curves, infinitely many of them have analytic rank in $\{2,4,6\}$.
\end{prop}
\begin{proof}
For the congruent number family we have $N_E=N_{E^1}=32$, but the root number is determined just by the class of $d$ modulo $8$ \cite{tunn}. In particular, it is $+1$ when $|\sqf(d)|\equiv1,2,3\pmod8$. Moreover, the homogeneous quartic of \S2.1 is $F(u,v)=uv(v+u)(v-u)$. Apply Definition \ref{def:avg} with $N=1$, $I=8\Z^++6$, $i:u\mapsto E^{F(u,1)}$ and $W\in C^1_c(\R)$ approximating $\mathbf{1}_{[0,1]}$, and let $\f=i(I)$. Then
\begin{align*}\avg r_{an}(\f) = \limsup_{T\to\infty}\frac{8}{T}\sum_{\substack{0\le u\le T \\ u\equiv_{8} 6}}r_{an}(E^{F(u,1)}),\end{align*} which we upper bound via \eqref{eq:expin}.

Let us treat the averages of the terms $\log N_{E^{F(u,1)}}$ and $U_i(E_{F(u,1)},X), \ i =1,2$ separately. Since $N_{E^d}=O(\sqf(d)^2)$, we get $\log N_{E^{F(u,1)}}\le 6\log u$. Therefore, when givided by $\log X,$ the average of their sum is \begin{align*}\limsup_{T\to\infty}\frac{8}{T}\sum_{\substack{u\le T \\ u\equiv_8 6}}\frac{\log N_{E^{F(u,1)}}}{\log X}\le6\frac{\log T}{\log X}.\end{align*}

The $U_2$ term contributes $\frac12+O_E(\frac{\log\log T}{\log X})$ thanks to Lemma \ref{lem:gold}.

Finally, we are left with the $U_1$ term. Notice that we can replace the Legendre symbol of $d=\sqf(F(u,1))$ that we get from Lemma \ref{lem:apa1} with that of $F(u,v)$: indeed, letting $F(u,1)=ds^2$, the two symbols agree except possibly at primes dividing $s$, and $\log|s|\ll\log(T)$. By the Hasse bound, the error for each pair $(u,v)$ is $\ll\sum_{p\mid s}(\log p)/\sqrt p$, and hence $\ll\sqrt{\log T}$. The only prime dividing $2N_E$ is $2$, so now the lemma holds for all primes we sum on; we have:
\begin{align}\begin{split}\lvert\sum_{\substack{u\le T \\ u\equiv_{8} 6}}U_1(E^{F(u,1)}&,X)\rvert\le\left\lvert\sum_{u}\sum_{\substack{5\le p\le X}}h_X(\log p)\frac{\log p}{p}a_p(E)\genfrac(){}{0}{F(u,1)}{p}\right\rvert \le \label{main}\end{split} \\ &\le \sum_{\substack{5\le p\le X}}h_X(\log p)\frac{\log p}{p}|a_p(E)|\left\lvert\sum_{u}\genfrac(){}{0}{F(u,1)}{p}\right\rvert.\notag\end{align}
Let us focus on the sum of Legendre symbols: the Chinese Remainder Theorem implies that
\begin{align}\left\lvert\sum_{u}\genfrac(){}{0}{F(u,1}{p}\right\rvert=\left\lvert\left\lfloor\frac{T}{8p}\right\rfloor a_p(E) + \sum_{x=p\left\lfloor\frac{T}{8p}\right\rfloor}^{\left\lfloor\frac{T-6}{8}\right\rfloor}\genfrac(){}{0}{F(8x+6,1)}{p}\right\rvert,\label{change}\end{align}since $y^2=F(u,1)$ is isomorphic to $E$.
As the range of values over which $x$ is running is shorter than $\frac{T}{8}-p(\frac{T}{8p}-1)=p$, Weil's bounds on sum of characters (see \cite[1.3]{cocpin}) combined with a classical generalization of the Pólya--Vinogradov inequality (see \cite[1.1]{bon}) imply, as Burgess states in \cite[Abstract]{burg}, that the modulus of the sum is $\ll \sqrt{p}\log p$, so we finally get:

\begin{align*}\vert\sum_{\substack{u\le T \\ u\equiv_{8} 6}}U_1(E_{F(u,1)},X)\vert\le T\sum_{5\le p\le X}\left(h_X(\log p)\frac{\log p}{8p^2}a_p(E)^2\right)+O\left(\sum_{5\le p\le X}p^{-\frac12}\log^2p|a_p(E)|\right).\end{align*}

The error term is $O(X\log X)$ by the Prime Number Theorem. For the main term we invoke \cite[Lemma 7.17]{bfkmms} (where $\lambda(p)=a_p/\sqrt p$), and partial summation on the series of $\frac{a_p^2\log p}{p^2}$ with the weight $h_X$ gives
\begin{align*}\frac{2}{\log X}\lvert\sum_{u}U_1(E^{F(u,1)},X)\rvert\le\frac{T}{8} + O(X).\end{align*}

Choosing $X\log X=T$ and letting $T\to\infty$ then gives
\begin{align*}\avg r_{an}(\f) \le \frac{1}{2}+6(1+\epsilon)+1\end{align*}for every $\epsilon>0$, so $\avg r_{an}(\f)<8$. The map $i$ is good since its fiber at $d$ consists of the $x$-coordinates of the integral points of $dy^2=x^3-x$, so Lemma \ref{lem:ineq} applies and gives infinitely many congruent number curves in $\f$ with analytic rank $<8$. As they all have positive analytic rank by Corollary \ref{posanal}, which is also even since $|\sqf(F(u,1)|\equiv2\pmod8$ for $u\equiv6\pmod8$, the claim follows.\end{proof}

\begin{proof}[Proof of Theorem \ref{odd}]
Let $E=E^1$ be the congruent number curve. Recall that, by the works of Gross--Zagier and Kolyvagin \cite{gz,kol}, the Birch--Swinnerton-Dyer conjecture is known for elliptic curves over $\Q$ with analytic rank smaller than 2. Therefore, if we have a curve with two independent nontorsion points, it cannot have analytic rank $\le1$. This allows us to apply the same idea of the proof of Proposition \ref{thmcn}: our goal is to find a positive-density subset $I\subset\Z^+$ and a function $i:I\to\f_E$ such that: \begin{itemize}
    \item the curves in $\g=i(I)$ have rank at least $2$ and root number $-1$;
    \item the map $i$ is good;
    \item $\avg r_{an}(\g)<19.$
\end{itemize}

The claim would then follow as in the proof of the proposition. In \cite{halhun}, Halbeisen and Hungerbühler prove that the curves $E_n$ with $$n=uv(u+2v)(2u+v)(u^2-v^2)(u^2+uv+v^2), \ u,v\in\Z$$ have rank at least $2$. Setting $v=1$, define $Q(u)=u(u-1)(u+1)(u+2)(2u+1)(u^2+u+1)$ and $Q_1(u)=\frac{Q(u)}{(u^2-1)}$. Notice that, for $u\equiv13\pmod{16}$, we have $Q_1(u)\equiv7\pmod8$ and $u^2-1\equiv8\pmod{32}$. 
Therefore, $\sqf(Q(u))\equiv7\cdot2\equiv6\pmod8$. Moreover, $Q$ has distinct roots, so the equation $dy^2=Q(u)$ defines a genus-$3$ curve, which has finitely many rational points by Faltings's theorem. Therefore, the set $I=\{u\in16Z^++13\}$ together with the map $i:u\mapsto E^{\mathrm{Q(u)}}$ satisfies the first two requirements.

The upper bound for the average analytic rank is then obtained in the same fashion as in the proof of Proposition \ref{thmcn} (in particular, notice how the application of the Chinese Remainder Theorem to sum over $\F_p$ is unchanged), with two differences: \begin{enumerate}
\item the expression for $n$ now has degree $7$ in our variable $u$, so the average of the $\log$ conductor terms for the same choice of $T=X^{1+\epsilon}$ gives us $14(1+\epsilon)$;
\item in the analog of the RHS of \eqref{change}, the first term is not $\left\lfloor\frac{T}{8p}\right\rfloor a_p(E)$, since the polynomial whose Legendre symbols we are summing is $Q$, and not the one defining a Weierstrass model for $E$. Nevertheless, the same bound \cite[1.3]{cocpin} mentioned above for the sum over the whole $\F_p$ shows that we can bound that term with $6\left\lfloor\frac{T}{8p}\right\rfloor\sqrt{p}$. As a consequence, we sum $6|a_p(E)|\sqrt{p}$ instead of $a_p(E)^2$, and the CM Sato--Tate distribution \cite{Sutherland2019} yields an extra factor of $\frac{12}{\pi}$.\end{enumerate}
Since the bound for $U_2$ is unchanged, we obtain:
\begin{align*}\avg r_{an}(\g)\le 14(1+\epsilon)+\frac{12}{\pi}+\frac{1}{2}<18.32,\end{align*}
which concludes the proof of Theorem \ref{odd}.

\end{proof}

\section{Proof of the main result}\label{sec:2dfourier}

We now turn our attention to Theorem \ref{mainthm}. Let us now work with a general base elliptic curve $E:y^2=f(x)$, with $f$ chosen to give an integral minimal model.

\subsection{Setup}

Let us return to the average of $U_1$. In \S\ref{cn} we exploited cancellation from the varying twisting parameter, but took absolute values before summing over primes. To improve on it, we need to use cancellation in $a_p(E)$ as well. Heath-Brown's result \cite[Lemma 6]{hb} gives such an estimate when $a_p(E)$ is multiplied by a quadratic character. Applying Poisson summation in one parameter to twists by $F(u,1)$, produces sums of the form $\sum_{x\bmod p}\bigl(\frac{f(x)}p\bigr)e_p(kx)$, to which his result does not apply. The idea to produce the desired character is to retain both parameters in $F(u,v)=u(v^3+au^2v+bu^3)$:
\begin{lem}\label{lem:2d-complete}
    Let $h,k$ be integers and let $p\nmid2N_E$ be prime. We have \begin{align}\label{eq:2d-complete}
 \sum_{u,v\bmod p}\genfrac(){}{0}{F(u,v)}p e_p(hu+kv)
 =p\genfrac(){}{0}{F(k,-h)}p+a_p(E)
   -p\,a_p(E)\mathbf1_{p\mid h,\ p\mid k}.
\end{align}
\end{lem}
\begin{proof}
    For $u\ne0$ modulo $p$, let $v=ux$. As $F(u,ux)=u^4f(x)$, the Legendre symbol does not depend on $u$, so the left-hand side of \eqref{eq:2d-complete} becomes
\[
 \sum_{x\bmod p}\genfrac(){}{0}{f(x)}p
   \sum_{u\ne0}e_p\bigl(u(h+kx)\bigr).
\]
The inner sum is $p-1$ when $h+kx=0$, and $-1$ otherwise. Since $\sum_x\genfrac(){}{0}{f(x)}p=-a_p(E)$ by the coprimality hypothesis, we obtain
\begin{align*}
\sum_{u,v\bmod p}\genfrac(){}{0}{F(u,v)}p e_p(hu+kv)
=p\sum_{\substack{x\bmod p\\h+kx=0}}
\genfrac(){}{0}{f(x)}p+a_p(E).
\end{align*}
If $p\nmid k$ the equation $h+kx=0$ has the only solution
$x=-h/k$, so the first term is
$p\genfrac(){}{0}{F(k,-h)}p$ as $F(k,-h)=k^4f(-h/k)$. If $p\mid k$ and $p\nmid h$ there are no solutions, and at the same time $\genfrac(){}{0}{F(k,-h)}p=0,$ so we again get \eqref{eq:2d-complete}. If $p$ divides both $h$ and $k$ there are $p$ solutions, giving exactly the extra term $-p\,a_p(E)\mathbf1_{p\mid h,\ p\mid k}$.
\end{proof}

Recall that we want to work with the twists $E^{F(u,v)}$ having root number $+1$ and for which the point of Lemma \ref{fam} is nontorsion. To select them, let $q=8N_E$ and set
\begin{equation}\label{eq:2d-local}
 I=\{(u,v)\in\Z^2:\ (v,q)=1,\ u\equiv_{q^3} q^2v,\ \omega(E)F(u,v)>0,\ \log\max\{|u|,|v|\}>C_E+\log\gcd(u,v)\}.
\end{equation}
For $(u,v)\in I$, we have $q^2\mid u$ and $u/q^2\equiv v\pmod q$, so $F(u,v)/q^2\equiv v^4\pmod q$. The squarefree part $d=\sqf(F(u,v))$ is therefore a square modulo $q$, forcing $d\equiv1\pmod8$ and $\chi_d(N_E)=1$. Modulo $q^3$, the pairs in $I$ are those of the form $(q^2t,t)$ with $(t,q)=1$. There are $\varphi(q^3)$ such pairs, so
\begin{align}\label{eq:2d-mass}
 \sum_{(u,v)\in I}W(u/T,v/T)
 =\frac{\varphi(q^3)}{q^6}T^2+O(T).\end{align}

Take a nonnegative $W\in C_c^\infty((-1,1)^2)$ supported in a small rectangle where $F$ has sign $\omega(E)$. Such a rectangle exists because, for fixed $v>0$ and sufficiently small $|u|$, $F(u,v)$ has the sign of $u$. The formula $\omega(E^d)=\omega(E)\operatorname{sgn}(d)\chi_d(N_E)$ \cite[p. 2]{hb} then gives root number $+1$ for every twist with nonzero weight. Let $\f=\{E^{F(u,v)}:(u,v)\in I\}.$ Corollary \ref{posanal} and Remark \ref{pos} then ensure that the curves in $\f$ have even analytic rank at least $2$.

\subsection{Key estimates}

We use $h(t)=\max(0,1-|t|)$, as in \eqref{eq:expin}. 
As in the proof of Proposition \ref{thmcn}, we can replace $\sqf(F(u,v))$ in the Legendre symbols obtained from Lemma \ref{lem:apa1} with just $F(u,v)$.
The finitely many primes dividing $q$ contribute a constant to $U_1$, so
\begin{align}\label{eq:2d-squarepart}
 &\sum_{(u,v)\in I}W(u/T,v/T)U_1(E^{F(u,v)},X)=\notag\\
 =-\sum_{\substack{5\le p\le X\\p\nmid q}}
 \frac{a_p(E)\log p}{p}&h_X(\log p)
 \sum_{(u,v)\in I}W(u/T,v/T)\genfrac(){}{0}{F(u,v)}p
 +O\bigl(T^2\sqrt{\log(2T)}\bigr).
\end{align}


As the congruence conditions defining $I$ are modulo $q^3$ and the Legendre symbol mod $p$, we split the inner sum in \eqref{eq:2d-squarepart} into residue classes modulo $q^3p$. Applying Poisson summation in both variables gives
\begin{align}
 &\sum_{(u,v)\in I}W(u/T,v/T)\genfrac(){}{0}{F(u,v)}p=\notag \\
 \quad=\frac{T^2}{q^6p^2}\sum_{h,k\in\Z}
 &\widehat W\left(\frac{Th}{q^3p},\frac{Tk}{q^3p}\right)
 \sum_{\substack{r,s\bmod q^3p\\
 r\equiv q^2s\ (\bmod q^3)\\(s,q)=1}}
 \genfrac(){}{0}{F(r,s)}p e_{q^3p}(hr+ks).\label{eq:poisson}
\end{align}
We now evaluate the sum over $r,s$. Since $p\nmid q$, each pair in the sum can be written uniquely as
\[
 r=pq^2t+q^3x,\qquad s=pt+q^3y\pmod{q^3p},
\]
where $t$ runs over the residues modulo $q^3$ coprime to $q$, and $x,y$ run over the residues modulo $p$. In particular,
\[
 e_{q^3p}(hr+ks)
 =e_{q^3}\bigl(t(q^2h+k)\bigr)e_p(hx+ky).
\]
Moreover, $F(r,s)\equiv q^{12}F(x,y)\pmod p$, so they have equal Legendre symbols. The inner sum in \eqref{eq:poisson} then factors as
\begin{align*}
 \left(\sum_{x,y\bmod p}
 \genfrac(){}{0}{F(x,y)}p e_p(hx+ky)\right)\sum_{\substack{t\bmod q^3\\(t,q)=1}}
 e_{q^3}\bigl(t(q^2h+k)\bigr).
\end{align*}

Now invoke Lemma \ref{lem:2d-complete} for the first factor and consider its three terms separately. The contribution of $a_p(E)$ is
\[
\frac{a_p(E)T^2}{q^6p^2}
\sum_{h,k\in\Z}
\widehat W\left(\frac{Th}{q^3p},\frac{Tk}{q^3p}\right)
\sum_{\substack{t\bmod q^3\\(t,q)=1}}
e_{q^3}\bigl(t(q^2h+k)\bigr)=a_p(E)\sum_{(u,v)\in I}W(pu/T,pv/T),
\]
where the equality follows from applying Poisson summation modulo $q^3$ to the weight $W(pu/T,pv/T)$.

For the term $-p\,a_p(E)\mathbf1_{p\mid h,\,p\mid k}$, set $h=ph'$ and $k=pk'$. Since multiplication by $p$ permutes the residues $t$ coprime to $q$, we have
\[
 \sum_{\substack{t\bmod q^3\\(t,q)=1}}
 e_{q^3}\bigl(pt(q^2h'+k')\bigr)
 =
 \sum_{\substack{t\bmod q^3\\(t,q)=1}}
 e_{q^3}\bigl(t(q^2h'+k')\bigr).
\]
Poisson summation modulo $q^3$ then shows that this term contributes
$-\frac{a_p(E)}p\sum_{(u,v)\in I}W(u/T,v/T)$. Together with the expression for the third term containing
$p\bigl(\frac{F(k,-h)}p\bigr)$, this gives:
\begin{align}\label{eq:2d-poisson}
 \sum_{(u,v)\in I}W(u/T,v/T)&
 \genfrac(){}{0}{F(u,v)}p=-\frac{a_p(E)}p
 \sum_{(u,v)\in I}W(u/T,v/T)
 +a_p(E)\sum_{(u,v)\in I}W(pu/T,pv/T)\\
 \qquad\quad+\frac{T^2}{q^6p}&
 \sum_{h,k\in\Z}
 \Biggl(\sum_{\substack{t\bmod q^3\\(t,q)=1}}
 e_{q^3}\bigl(t(q^2h+k)\bigr)\Biggr)
 \genfrac(){}{0}{F(k,-h)}p
 \widehat W\left(\frac{Th}{q^3p},\frac{Tk}{q^3p}\right)\notag.
\end{align}

Let us treat the three terms in order. For the first one, notice that Lemma \ref{lem:gold} gives $U_2(E,X^2)=\tfrac12\log X+O(1)$. As the coefficient $a_{p^2}(E)$ used in $U_2$ is $a_p(E)^2-2p$ at good primes, and $h_{X^2}(2\log p)=h_X(\log p)$, we have
\begin{align}\label{eq:2d-mainconstant}
 \sum_{\substack{5\le p\le X\\p\nmid q}}
 \frac{a_p(E)^2\log p}{p^2}h_X(\log p)=2\sum_{p\le X}\frac{\log p}{p}h_X(\log p)-U_2(E,X^2)+O(1)=\frac12\log X+O(1),
\end{align}
by the Prime Number Theorem. Since the sum of the weights is $O(T^2)$ by \eqref{eq:2d-mass}, the total contribution is
\[
\frac12\log X\sum_{(u,v)\in I}W(u/T,v/T)+O(T^2).
\]

Similarly, the second term on the right-hand side of \eqref{eq:2d-poisson} gives
\[
-\sum_{\substack{5\le p\le X\\p\nmid q}}
\frac{a_p(E)^2\log p}{p}h_X(\log p)
\sum_{(u,v)\in I}W(pu/T,pv/T).
\]
Since $W$ is supported in $(-1,1)^2$ and $v\ne0$ for $(u,v)\in I$, the inner sum vanishes for $p\ge T$. For $p<T$, it is $O((T/p)^2)$, because only pairs with $|u|,|v|<T/p$ can contribute. As $a_p(E)^2\le4p$ and $|h_X|\le1$, the absolute value of this contribution is
\[
\ll T^2\sum_p\frac{\log p}{p^2}\ll T^2.
\]

For the last term, we interchange the sum over primes and that over $h,k$. \cite[Lemma 6]{hb} gives, for every nonzero integer $n$ and every $\epsilon>0$,
\begin{equation}\label{eq:2d-hb}
 \sum_{5\le p\le X}\frac{a_p(E)\log p}{p}\genfrac(){}{0}{n}p
 \ll_{\epsilon}X^\epsilon|n|^\epsilon.
\end{equation}
Indeed, the lemma is actually stated for the primitive character $\chi_{\sqf(n)}$, but the difference is supported on primes dividing the square factor of $n$, whose contribution is $O_\epsilon(|n|^\epsilon)$ by the Hasse bound. We split the prime sum into dyadic intervals
$P<p\le\min(2P,X)$, where $P$ runs through powers of $2$.
For fixed $h,k$ with $F(k,-h)\ne0$, partial summation
using \eqref{eq:2d-hb} gives
\begin{align}\label{eq:2d-frequency}
&\left|\sum_{\substack{P<p\le\min(2P,X)\\p\nmid q}}
 \frac{a_p(E)\log p}{p^2}\genfrac(){}{0}{F(k,-h)}p
 h_X(\log p)\widehat W\left(\frac{Th}{q^3p},\frac{Tk}{q^3p}\right)\right|\ll_{\epsilon,A}\notag\\
& \ll_{\epsilon,A}P^{-1+\epsilon}|F(k,-h)|^\epsilon
 \left(1+\frac{T(|h|+|k|)}P\right)^{-A}
\end{align}
for any $A>0$. 
This is where the gain with respect to the one-variable argument is realized: retaining the full Legendre symbol $\genfrac(){}{0}{F(k,-h)}p$ allows us to apply Heath-Brown's Lemma to the quadratic character obtained by fixing one of the variables, which saves $P^{1/2-\epsilon}$ in this prime sum.

For $P\ge T$ and assuming $A>2+4\epsilon$, summing the right-hand side of \eqref{eq:2d-frequency} over $h,k$ and including the factor $T^2$ in \eqref{eq:2d-poisson} gives a bound of order
\begin{equation}\label{eq:2d-block}
 T^2P^{-1+\epsilon}(P/T)^{2+4\epsilon}
 =P^{1+5\epsilon}T^{-4\epsilon}\le P^{1+5\epsilon},
\end{equation}
since $|F(k,-h)|\ll(|h|+|k|)^4$. For $P<T$, the same sum is $O(T^2P^{-1+\epsilon}(P/T)^A)=O(P^{1+\epsilon})$. Summing over dyadic intervals gives the upper bound $O_\epsilon(X^{1+\epsilon})$ for the last term's contribution. 

Combining the three terms, the error from removing square factors, and \eqref{eq:2d-mass}, we obtain:
\begin{align}\label{eq:2d-u1average}
 \frac{\displaystyle\sum_{(u,v)\in I}W(u/T,v/T)
 U_1(E^{\sqf(F(u,v))},X)}
 {\displaystyle\sum_{(u,v)\in I}W(u/T,v/T)}
 =\frac12\log X+O(\sqrt{\log T})
 +O_{\epsilon}\left(\frac{X^{1+\epsilon}}{T^2}\right).
\end{align}

\subsection{Proof of Theorem \ref{mainthm}}
Recall that every twist in $\f$ has even analytic rank at least $2$. Also, by \eqref{eq:expfor}, for any $X\ge2$ the quantity $(\log X)\mathrm{avg}_Wr_{\mathrm{an}}(\f)$ is upper bounded by the average log conductor of $\f$, plus twice its average $U_1(X)$ and $U_2(X)$.
Setting $X=T^{2-\delta}$, $0<\delta<2$, and choosing $\epsilon>0$ such that $(2-\delta)(1+\epsilon)<2$, makes \eqref{eq:2d-u1average} contribute $\log X(1+o(1))$ to this. The contribution of $2U_2(X)$ is again at most $\log X(\frac12+o(1))$ by Lemma \ref{lem:gold}. Finally, the conductor bound $N_{E^d}\ll_Ed^2$, together with $|F(u,v)|\ll T^4$, gives:
\begin{equation}\label{eq:2d-rankaverage}
 \mathrm{avg}_Wr_{\mathrm{an}}(\f)
 \le\frac8{2-\delta}+1+\frac12.
\end{equation}
Therefore, taking $0<\delta<2/9$ in \eqref{eq:2d-rankaverage}, we obtain $\mathrm{avg}_Wr_{\mathrm{an}}(\f)<6.$

In order to apply Lemma \ref{lem:ineq}, it remains to check that our $i$ is good. Fix a nonzero squarefree integer $d$. Every pair with $\sqf(F(u,v))=d$ satisfies $\genfrac(){}{0}{F(u,v)}p\in\{0,\genfrac(){}{0}{d}p\}$ for every prime $p\nmid2dN_E$. There are at most $4p$ pairs modulo $p$ with $F(u,v)=0$, and \eqref{eq:2d-complete} at $h=k=0$ gives $\sum_{u,v\bmod p}\genfrac(){}{0}{F(u,v)}p=-(p-1)a_p(E)$. By the Hasse bound, it follows that
\begin{align*}
 \frac1{p^2}\#\left\{(u,v)\in\F_p^2:\genfrac(){}{0}{F(u,v)}p\in\left\{0,\genfrac(){}{0}{d}p\right\}\right\}
 =\frac12+O(p^{-1/2}).
\end{align*}
For $p$ sufficiently large this is less than $3/4$. Applying the Chinese Remainder Theorem to any $j$ distinct such primes gives an upper bound of $(3/4)^j$ for the upper density of $\{(u,v)\in I:\sqf(F(u,v))=d\}$, so its density is $0$. Lemma \ref{lem:ineq} then applies, giving infinitely many curves in $i(I)$ with analytic rank smaller than $6$, and we are done.

\printbibliography
\end{document}